# 2-Clean Rings [*]

Z. Wang and J.L. Chen

**Abstract.** A ring $R$ is said to be $n$-clean if every element can be written as a sum of an idempotent and $n$ units. The class of these rings contains clean ring and $n$-good rings in which each element is a sum of $n$ units. In this paper, we show that for any ring $R$, the endomorphism ring of a free $R$-module of rank at least 2 is 2-clean and that the ring $B(R)$ of all $\omega \times \omega$ row and column-finite matrices over any ring $R$ is 2-clean. Finally, the group ring $RC_n$ is considered where $R$ is a local ring.

## 1. Introduction

The question of when the automorphism group of a module additively generates its endomorphism ring has been of interest for many years. A ring is called $n$-good [12] if every element is a sum of $n$ units. In 1953 Wolfson [14] and in 1954 Zelinsky [17] showed, independently, that every element of the ring of all linear transformations of a vector space over a division ring of characteristic not 2 is 2-good. In 1985 Goldsmith [4] proved that the endomorphism ring of a complete module over a complete discrete valuation ring is 2-good. In [13] Wans considered free $R$-modules where $R$ is a $PID$, and showed that if the rank of $M$ is finite and greater than 1, then $End_R(M)$ is 2-good. Meehan [8] further showed that the endomorphism ring of a free $R$-module of rank at least 2 is 2-good where $R$ is a PID. Moreover, the above question is considered by many authors on abelian groups (see [2],[7],[8]) and on general ring with an identity (see [3],[6],[11]).

[*]This work was supported by the Foundation for Excellent Doctoral Dissertation of Southeast University (YBJJ0507), the National Natural Science Foundation of China (No.10571026) and the Natural Science Foundation of Jiangsu Province (No.BK2005207).

In 1977 Nicholson [10] introduced the concept of a clean ring (1-clean) which contains unit-regular rings and semiperfect rings, and showed that every clean ring must be exchange. Camillo and Yu [1] further proved that a clean ring with 2 invertible is 2-good. Recently, Xiao and Tong [16] called a ring $R$ $n$-clean if every element of $R$ is the sum of an idempotent and $n$ units. The class of these rings contains clean rings and $n$-good rings. In 1974 Henriksen [6] found that for any ring $R$ and $n > 1$, the matrix ring $M_n(R)$ is 3-good. Moreover, Vámos [12] proved that for any ring $R$, the endomorphism ring of a free $R$-module of rank at least 2 is 3-good. Motivated by the result of Henriksen and Vámos, we conjecture that for any ring $R$, the endomorphism ring of a free $R$-module of rank at least 2 is 2-clean.

In this paper, we answer the question in the positive. In fact, we proved that for any ring $R$, the endomorphism ring of a free $R$-module of rank at least 2 is 2-clean. It is also proved that the ring $B(R)$ of all $\omega \times \omega$ row and column-finite matrices over any ring $R$ is 2-clean. Finally, the group ring $RC_n$ is considered where $R$ is a local ring.

Throughout this paper, rings are associative with identity and modules are unitary. $J(R)$ and $U(R)$ denote the Jacobson radical and the group of units of $R$, respectively.

## 2. BASIC PROPERTIES OF $n$-CLEAN RINGS

An element of a ring is called $n$-clean if it can be written as the sum of an idempotent and $n$ units. A ring is called $n$-clean if each of its elements is $n$-clean. In this section, some properties of $n$-clean rings are given.

**Proposition 1.** *Let $R$ be a ring and let $a \in R$. Then the following statements hold:*

*(1) if $a$ is $n$-clean then it is also $l$-clean for all $n \leq l$.*

*(2) every $n$-good ring is $n$-clean; if $R$ is $n$-clean with $2 \in U(R)$ then it is $(n+1)$-good.*

**Proof.** (1) We only need to prove that $a$ is $n+1$-clean. Let $a \in R$ be $n$-clean: $a = e + u_1 + u_2 + \cdots + u_n$ where $e^2 = e \in R$ and $u_1, u_2, \cdots, u_n \in U(R)$. Note that $e = (1-e) + (2e-1)$, thus we have $a = (1-e) + (2e-1) + u_1 + \cdots + u_n$ where $2e - 1 \in U(R)$.

(2) It is clear that every $n$-good ring is $n$-clean. The second statement is due to Xiao and Tong (see [16]). □

Let $S(R)$ be the nonempty set of all proper ideal of $R$ generated by central idempotents. An ideal $P \in S(R)$ is called a Pierce ideal of $R$ if $P$ is a maximal (with respect to inclusion) element of the set $S(R)$. If $P$ is a Pierce ideal of $R$, then the factor ring $R/P$ is called a Pierce stalk of $R$. The next result shows that the $n$-clean property needs to be checked only by for indecomposable rings or Pierce stalks.

**Proposition 2.** *Let $R$ be a ring. Then the following are equivalent:*

(1) *$R$ is $n$-clean.*

(2) *every factor ring of $R$ is $n$-clean.*

(3) *every indecomposable factor ring of $R$ is $n$-clean.*

(4) *every Pierce stalk of $R$ is $n$-clean.*

**Proof.** (1) $\Rightarrow$ (2), (2) $\Rightarrow$ (3) and (2) $\Rightarrow$ (4) are directly verified.

(3) $\Rightarrow$ (1). Suppose that (3) holds and $R$ is not $n$-clean, then there is an element $a \in R$ which is not $n$-clean. Now let $\mathcal{S}$ be the set of all proper ideals $I$ of $R$ such that $\overline{a}$ is not $n$-clean in $R/I$. Clearly, $0 \in \mathcal{S}$ and the set $\mathcal{S}$ is not empty. Define a partial ordering on $\mathcal{S}$ by $"\subseteq"$. If $\{I_\alpha : \alpha \in \Lambda\}$ is a chain in $\mathcal{S}$, let $I = \cup_{\alpha \in \Lambda} I_\alpha$. We will show that $\overline{a}$ is not $n$-clean in $R/I$. Suppose that $\overline{a}$ is $n$-clean in $R/I$. Then there exist $\overline{u}_1, \overline{u}_2, \cdots, \overline{u}_n \in U(R/I)$ (with inverses $\overline{v}_1, \overline{v}_2, \cdots, \overline{v}_n$, respectively) and $\overline{e}^2 = \overline{e} \in R/I$ such that $\overline{a} = \overline{e} + \overline{u}_1 + \overline{u}_2 + \cdots + \overline{u}_n$. Note that $e^2 - e \in \cup_{\alpha \in \Lambda} I_\alpha$ and $u_i v_i - 1$, $v_i u_i - 1 \in \cup_{\alpha \in \Lambda} I_\alpha$, so $e^2 - e \in I_{\alpha_0}$, $u_i v_i - 1 \in I_{\alpha_i}$ and $v_i u_i - 1 \in I_{\alpha_i'}$ for $\alpha_0, \alpha_i, \alpha_i' \in \Lambda$. Because $\{I_\alpha : \alpha \in \Lambda\}$ is a chain in $\mathcal{S}$, there is a maximal $I_s$ in the set $\{I_{\alpha_0}, I_{\alpha_1}, \cdots, I_{\alpha_n}, I_{\alpha_1'}, I_{\alpha_1'}, \cdots, I_{\alpha_n'}\}$ such that $I_{\alpha_0}, I_{\alpha_i}, I_{\alpha_i'} \subseteq I_s$. That is, $\overline{a}$ is $n$-clean in $R/I_s$, a contradiction. This implies that $I \in \mathcal{S}$ is a upper bound of the chain. Because $\mathcal{S}$ is an inductive set and, by Zorn's Lemma, $\mathcal{S}$ has a maximal element $I_0$. By (3) $R/I_0$ is decomposable as a ring. Write $R/I_0 \cong R/I_1 \oplus R/I_2$ where both the ideals $I_1$, $I_2$ strictly contain $I_0$ and so by the choice of $I_0$, $\overline{a}$ is $n$-clean in $R/I_1$ and $R/I_2$. But then $\overline{a}$ is $n$-clean in $R/I_0$, a contradiction.

(4) $\Rightarrow$ (1). Let $\mathcal{S}$ be the set of all proper ideals $I$ of $R$ such that $I$ is generated by central idempotents and the ring $R/I$ is not $n$-clean. Assume that $R$ is not $n$-clean. Then $0 \in \mathcal{S}$ and the set $\mathcal{S}$ is not empty. It is directly verified as above that the union of every ascending chain of ideals from $\mathcal{S}$ belongs to $\mathcal{S}$. By Zorn's Lemma, the set $\mathcal{S}$ contains a maximal element $P$. By condition (4), it is sufficient to prove that $P$ is a Pierce ideal. Assume that

contrary. By the definition of the Pierce ideal, there is a central idempotent $e$ of $R$ such that $P+eR$ and $P+(1-e)R$ are proper ideals of $R$ which properly contain the ideal $P$. Since ideals $P+eR$ and $P+(1-e)R$ do not belong to $\mathcal{S}$ and are generated by central idempotents, $R/(P+eR)$ and $R/(P+(1-e)R)$ are $n$-clean. Note that $R/P \cong (R/(P+eR)) \times (R/(P+(1-e)R))$, it can be verified that $R$ is $n$-clean. □

## 3. MATRIX RINGS AND ENDOMORPHISM RINGS OF FREE MODULES

In this section, we will consider the 2-cleaness of the endomorphism ring of a free $R$-module of rank at least 2. First we give the following simple and interesting decomposition.

**Lemma 3.** *Over any ring, the $2 \times 2$ and $3 \times 3$ matrices are 2-clean.*

**Proof.** Let $R$ be a ring and let $A = \begin{pmatrix} a_{11} & a_{12} \\ a_{21} & a_{22} \end{pmatrix} \in M_2(R)$. Put $E = \begin{pmatrix} a_{11} - 1 & 2 - a_{11} \\ a_{11} - 1 & 2 - a_{11} \end{pmatrix}$. It is checked easily that then $E^2 = E$. Thus we have

$$A - E = \begin{pmatrix} 1 & a_{12} + a_{11} - 2 \\ a_{21} - a_{11} + 1 & a_{22} + a_{11} - 2 \end{pmatrix}.$$

Observing the above matrix, and then there exist invertible matrices

$$P = \begin{pmatrix} 1 & 0 \\ a_{11} - a_{21} - 1 & 1 \end{pmatrix} \quad \text{and} \quad Q = \begin{pmatrix} 1 & 2 - a_{11} - a_{12} \\ 0 & 1 \end{pmatrix}$$

such that

$$P(A - E)Q = \begin{pmatrix} 1 & 0 \\ 0 & c \end{pmatrix} = \begin{pmatrix} 1 & 1 \\ 1 & 0 \end{pmatrix} + \begin{pmatrix} 0 & -1 \\ -1 & c \end{pmatrix},$$

where $c = a_{11}^2 + a_{11}a_{12} - a_{21}a_{12} - a_{21}a_{11} - 2a_{11} + 2a_{21} - a_{12} + a_{22}$. This shows that $A = P^{-1} \begin{pmatrix} 1 & 1 \\ 1 & 0 \end{pmatrix} Q^{-1} + P^{-1} \begin{pmatrix} 0 & -1 \\ -1 & c \end{pmatrix} Q^{-1} + E$ is 2-clean.

Now let $B = \begin{pmatrix} b_{11} & b_{12} & b_{13} \\ b_{21} & b_{22} & b_{23} \\ b_{31} & b_{32} & b_{33} \end{pmatrix}$ be a $3 \times 3$ matrix over R. We first construct an idempotent in order to show 2-cleaness of $B$. Set

$$F = \begin{pmatrix} b_{11} - 1 & b_{22} - 1 & 3 - b_{11} - b_{22} \\ b_{11} - 1 & b_{22} - 1 & 3 - b_{11} - b_{22} \\ b_{11} - 1 & b_{22} - 1 & 3 - b_{11} - b_{22} \end{pmatrix}.$$

It is directly verified that $F^2 = F$. Thus

$$B - F = \begin{pmatrix} 1 & b_{12} - b_{22} + 1 & b_{13} + b_{11} + b_{22} - 3 \\ b_{21} - b_{11} + 1 & 1 & b_{23} + b_{11} + b_{22} - 3 \\ b_{31} - b_{11} + 1 & b_{32} - b_{22} + 1 & b_{33} + b_{11} + b_{22} - 3 \end{pmatrix}.$$

We only need to show that $B - F$ is 2-good. Observing the above matrix, and then there exist invertible matrices

$$T = \begin{pmatrix} 1 & 0 & 0 \\ 0 & 1 & 0 \\ b_{11} - b_{31} - 1 & 0 & 1 \end{pmatrix}, V = \begin{pmatrix} 1 & b_{22} - b_{12} - 1 & 0 \\ 0 & 1 & 0 \\ 0 & 0 & 1 \end{pmatrix}$$

and $W = \begin{pmatrix} 1 & 0 & 0 \\ 0 & 1 & 3 - b_{23} - b_{11} - b_{22} \\ 0 & 0 & 1 \end{pmatrix}$ such that

$$VT(B - F)W = \begin{pmatrix} * & 0 & * \\ * & 1 & 0 \\ 0 & * & * \end{pmatrix} = \begin{pmatrix} 0 & 1 & * \\ 0 & 0 & 1 \\ 1 & * & * \end{pmatrix} + \begin{pmatrix} * & -1 & 0 \\ * & 1 & -1 \\ -1 & 0 & 0 \end{pmatrix}.$$

Consider the two matrices $U_1$, $U_2$ occurring in the decomposition above of $VT(B - F)W$. It is straightforward to verify that the two matrices are invertible in $M_3(R)$. Thus we obtain immediately a 2-clean expression of $B$, i.e.,

$$B = T^{-1}V^{-1}U_1W^{-1} + T^{-1}V^{-1}U_2W^{-1} + F.$$

This completes the proof. □

**Remark 4.** (1). For the matrix ring $M_n(R)$, it is customary to write $GL_n(R)$ for $U(M_n(R))$. An elementary matrix is the result of an elementary row operation performed on the identity matrix. We denote by $E_n(R)$ the subgroup of $GL_n(R)$ generated by the elementary matrices, permutation matrices and -1. Observing the decompositions of the $2 \times 2$ and $3 \times 3$ matrices above, we see that, these matrices can be written as the sum of an idempotent matrix and two elements of $E_n(R)$.

(2). For any ring $R$, $R$ can be embedded in the $2 \times 2$ matrix ring $M_2(R)$. That is, all rings can be embedded in a 2-clean ring by Lemma 3.

(3). We know that 2-clean rings contain clean rings and 2-good rings. However, the converse is not true. For example, the matrix ring $M_2(\mathbb{Z})$ is not clean since $\mathbb{Z}$ is not a exchange ring, and the matrix ring $M_2(\mathbb{Z}[x])$ is not 2-good (see [12, Proposition 8]).

(4). It is well known that for a clean ring $R$, idempotents can be lifted modulo $J(R)$. However, a 2-clean ring has not this property in general. Let $R = \mathbb{Z}_{(2)} \cap \mathbb{Z}_{(3)} = \{m/n \in \mathbb{Q} : m, n \in \mathbb{Z}, 2 \nmid n \text{ and } 3 \nmid n\}$ and set $S = M_2(R)$. Then $J(S) = J(M_2(R)) = M_2(J(R)) = M_2(6R)$. Let $F = \begin{pmatrix} 3 & 0 \\ 6 & 3 \end{pmatrix}$. Then $F^2 - F \in J(S)$, but there is no idempotent $E$ of $S$ such that $F - E \in J(S)$ since non-trivial idempotents of $S$ are only of form $\begin{pmatrix} a & b \\ c & 1-a \end{pmatrix}$ where $bc = a - a^2$ for $a, b, c \in R$. Thus $S$ is 2-clean by Lemma 3 but there exists an idempotent which can not be lifted modulo $J(S)$.

**Lemma 5.** *Let $R$ be a ring, $m, n \geq 1$ and $k \geq 2$. If the matrix rings $M_n(R)$ and $M_m(R)$ are both k-clean, then so is the matrix ring $M_{n+m}(R)$.*

**Proof.** Let $A \in M_{n+m}(R)$ be a typical $(n+m) \times (n+m)$ matrix which we will write in the block decomposition form

$$A = \begin{pmatrix} A_{11} & A_{12} \\ A_{21} & A_{22} \end{pmatrix},$$

where $A_{11} \in M_n(R)$, $A_{22} \in M_m(R)$ and $A_{12}, A_{22}$ are appropriately sized rectangular matrices. By hypothesis, there exist invertible $n \times n$, $m \times m$ matrices $U_1, U_2, \cdots, U_k$ and $V_1, V_2, \cdots, V_k$, and idempotent matrices $E_1, E_2$ such that $A_{11} = U_1 + U_2 + \cdots + U_k + E_1$ and $A_{22} = V_1 + V_2 + \cdots + V_k + E_2$. Thus the decomposition

$$\begin{pmatrix} A_{11} & A_{12} \\ A_{21} & A_{22} \end{pmatrix} = \begin{pmatrix} U_1 & A_{12} \\ O & V_1 \end{pmatrix} + \begin{pmatrix} U_2 & O \\ A_{21} & V_2 \end{pmatrix} + \cdots + \begin{pmatrix} U_k & O \\ O & V_k \end{pmatrix} + \begin{pmatrix} E_1 & O \\ O & E_2 \end{pmatrix}$$

shows that $A$ is $k$−clean. □

**Corollary 6.** *Let $k \geq 1$. If $R$ is a k-clean ring, then so the matrix ring $M_n(R)$ for any positive integer $n$.*

**Proof.** For $k = 1$, it follows from [5, Corollary 1]. Assume that $k \geq 2$, it is clear by induction and by Lemma 5. □

**Theorem 7.** *Let $R$ be a ring and let the free $R$-module $F$ be (isomorphic to) the direct sum of $\alpha \geq 2$ copies of $R$ where $\alpha$ is a cardinal number. Then the ring of endomorphisms $E$ of $F$ is 2-clean.*

**Proof.** Assume first that $\alpha \geq 2$ is finite so $E \cong M_\alpha(R)$. Then $E$ is 2-clean for $\alpha = 2, 3$ by Lemma 3 and the values of $\alpha < \omega$ for which $E$ is 2-clean are closed under addition by Lemma 5. So $E$ is 2-clean for all finite $\alpha$.

Assume now that $\alpha$ is infinite. Then $E \cong M_2(E)$ follows from $F \cong F \oplus F$, and so $E$ is 2-clean by Lemma 3. □

## 4. ROW AND COLUMN-FINITE MATRIX RINGS

Let $B(R)$ be the ring of all $\omega \times \omega$ row and column-finite matrices over a ring $R$. Fix a free $R$-module $F = \bigoplus_{i=1}^{\infty} f_i R$ on countably many generators, and for each $k \in \mathbb{N}$ let $F_k = \bigoplus_{i=k}^{\infty} f_i R$. A moment's reflection, using the standard correspondence between $R$-endomorphisms of $F_R$ and $\omega \times \omega$ column-finite matrices over $R$ relative to the basis $\{f_i\}_{i=1}^{\infty}$, confirms that

$$B(R) \cong \{\phi \in End_R(F) : \text{for each } k \in \mathbb{N}, \exists\, m \in \mathbb{N} \text{ with } \phi(F_m) \subseteq F_k\}.$$

Hence we identify $B(R)$ with this ring of transformations. Next we will consider the 2-cleanness of $B(R)$. The proof of the following result is a modification of that in [8, Theorem 3.5].

**Theorem 8.** *Let $R$ is ring. Then the row and column-finite matrix ring $B(R)$ is 2-clean.*

**Proof.** Note that $B(R) \cong B(M_2(R))$, so we may assume that $R$ is 2-clean by Lemma 3. Let $\phi \in B(R)$. Recall that $\varphi$ is defined by

(a) $\alpha$-endomorphism if $\varphi(f_i R) \subseteq \bigoplus_{k>i} f_i R$ for all $i < \omega$;
(b) $\beta$-endomorphism if $\varphi(f_i R) \subseteq \bigoplus_{k=1}^{i-1} f_i R$ for all $i < \omega$;
(c) $d$-endomorphism if $\varphi(f_i R) \subseteq f_i R$ for all $i < \omega$.

Then $\phi$ can obviously be expressed as

$$\phi = \eta + \rho + \delta,$$

where $\eta$ is an $\alpha$-endomorphism, $\rho$ is a $\beta$-endomorphism and $\delta$ is a $d$-endomorphism. Since $\phi \in B(R)$, for each $k \in \mathbb{N}$, there exists $m \in \mathbb{N}$ such that $\phi(F_m) \subseteq F_k$. By the definitions of $\eta$, $\rho$ and $\delta$, we check easily that $\eta(F_m) \subseteq F_k$, $\rho(F_m) \subseteq F_k$ and $\delta(F_m) \subseteq F_k$. For the $\alpha$-endomorphism $\eta$, by [8, Proposition 3.2], there exists a strictly ascending sequence of integers $0 < r_0 < r_1 < r_2 < \cdots$ such that $\eta(f_i R) \subseteq \bigoplus_{k=i+1}^{s+2-1} f_k R$ for all $r_s \leq i < r_{s+1}$. Using this sequence we define endomorphisms $\eta_1$, $\eta_2$ of $F$ as follows

$$\eta_1 f_i = \begin{cases} \eta f_i & \text{for } r_{2t} \leq i < r_{2t+1}; \\ 0 & \text{for } r_{2t+1} \leq i < r_{2t+2}, \end{cases}$$

and

$$\eta_2 f_i = \begin{cases} 0 & \text{for } r_{2t} \leq i < r_{2t+1}; \\ \eta f_i & \text{for } r_{2t+1} \leq i < r_{2t+2}. \end{cases}$$

Clearly, $\eta_1$ and $\eta_2$ are $\alpha$-endomorphisms of $F$ with $\eta = \eta_1 + \eta_2$, and for each $k \in \mathbb{N}$, there exists $m \in \mathbb{N}$ such that $\eta_1(F_m) \subseteq F_k$ and $\eta_2(F_m) \subseteq F_k$. By [8,

Lemma 3.4], we have that $\eta_1$, $\eta_2$ are both locally nilpotent. Next we decompose the $\beta$-endomorphism $\rho$. For each $i < \omega$, we have

$$\rho f_i = \sum_{k<i} f_k r_{ik} = \sum_{\substack{k<i \\ k \in I_1}} f_k r_{ik} + \sum_{\substack{k<i \\ k \in I_2}} f_k r_{ik},$$

where $I_1 = \bigcup_{t<\omega}\{k \mid r_{2t} \leq k < r_{2t+1}\}$ and $I_2 = \bigcup_{t<\omega}\{k \mid r_{2t+1} \leq k < r_{2t+2}\}$. We define $\rho_1$, $\rho_2$ correspondingly, i.e.,

$$\rho_1 f_i = \sum_{\substack{k<i \\ k \in I_1}} f_k r_{ik} \quad \text{and} \quad \rho_2 f_i = \sum_{\substack{k<i \\ k \in I_2}} f_k r_{ik}.$$

Clearly, $\rho = \rho_1 + \rho_2$ and $\rho_1$, $\rho_2$ are both locally nilpotent. We check easily that for each $k \in \mathbb{N}$, there exists $m \in \mathbb{N}$ such that $\rho_1(F_m) \subseteq F_k$ and $\rho_2(F_m) \subseteq F_k$. Note that $\rho_1 \eta_2 = 0 = \rho_2 \eta_1$ by definitions of $\eta_1$, $\eta_2$, $\rho_1$, $\rho_2$, so $\eta_1 + \rho_2$ and $\eta_2 + \rho_1$ are also locally nilpotent. Now we consider the $d$-endomorphism $\delta$. For each $i < \omega$, there exists an element $r_i$ of $R$ such that $\delta f_i = f_i r_i$. Since $R$ is 2-clean, there are $e_i^2 = e_i \in R$ and units $u_{i1}$, $u_{i2}$ of $R$ such that

$$\delta f_i = f_i u_{i1} + f_i u_{i2} + f_i e_i.$$

defining $\delta_e f_i = f_i e_i$ and $\delta_j f_i = f_i u_{ij}$ ($i < \omega$, $j = 1, 2$). So $\delta = \delta_1 + \delta_2 + \delta_e$ and $\delta_1$, $\delta_2$, $\delta_e$ are $d$-endomorphisms of $F$. Note that for each $k \in \mathbb{N}$, set $m = k$, we get $\delta_1(F_m) \subseteq F_k$, $\delta_2(F_m) \subseteq F_k$ and $\delta_e(F_m) \subseteq F_k$. Thus we consider the decomposition of $\phi$

$$\begin{aligned}
\phi &= \eta + \rho + \delta \\
&= \eta_1 + \eta_2 + \rho_1 + \rho_2 + \delta_1 + \delta_2 + \delta_e \\
&= (\eta_1 + \rho_2 + \delta_1) + (\eta_2 + \rho_1 + \delta_2) + \delta_e \\
&= \delta_1(\delta_1^{-1}(\eta_1 + \rho_2) + 1) + \delta_2(\delta_2^{-1}(\eta_2 + \rho_1) + 1) + \delta_e.
\end{aligned}$$

Note that $\delta_1^{-1}(\eta_1 + \rho_2)$ is locally nilpotent since $\delta_1^{-1}$ is $d$-endomorphism and $\eta_1 + \rho_2$ is locally nilpotent, and so $\delta_1^{-1}(\eta_1 + \rho_2) + 1$ is an automorphism of $F$. Hence $\delta_1(\delta_1^{-1}(\eta_1 + \rho_2) + 1)$ is also an automorphism of $F$. Similarly, $\delta_2(\delta_2^{-1}(\eta_2 + \rho_1) + 1)$ is an automorphism of $F$. Clearly, by the definitions of $\delta_e$, $\delta_e$ is idempotent endomorphism of $F$. It is checked easily that $\eta_1 + \rho_2 + \delta_1$, $\eta_2 + \rho_1 + \delta_2$, $\delta_e \in B(R)$ since $B(R)$ is a ring. Thus we complete the proof. $\square$

**Remark 9.** From the proof of Theorem 8, we may consider row and column-finite matrix rings over a 2-good ring similarly. In fact, we obtain that if $R$ is 2-good then so is the row and column-finite matrix ring $B(R)$, and that for any ring $R$ the row and column-finite matrix ring $B(R)$ is 3-good.

## 5. 2-CLEAN GROUP RINGS

Given a group $G$ and a ring $R$, denote the group ring by $RG$. In this section, we consider the group ring $RC_n$ where $R$ is a local ring and $C_n$ is a cyclic group of order $n$. Some results of Xiao and Tong [16] are extended.

**Theorem 10.** *Let $R$ be a local ring with $\overline{R} = R/J(R)$ and let $C_n$ be a cyclic group of order $n$. If $char\overline{R} \neq 2$, then $RC_n$ is 2-good.*

**Proof.** If $char\overline{R} = 0$ or $(char\overline{R},\ n) = 1$, then $\overline{n}$ and $\overline{2}$ are invertible in $\overline{R}$. Note that $\overline{R}$ is a division ring, then $\overline{R}C_n$ is semisimple from $n \cdot \overline{1} = \overline{n} \in U(\overline{R})$, and so $\overline{R}C_n$ is clean. This implies that $\overline{R}C_n$ is 2-good by [1, Proposition 10]. We know that if $G$ is locally finite then $J(R)G \subseteq J(RG)$ by [15]. Clearly, $J(R)C_n \subseteq J(RC_n)$, and then $\overline{R}C_n \cong RC_n/J(R)C_n \twoheadrightarrow RC_n/J(RC_n)$. So the factor ring $RC_n/J(RC_n)$ is 2-good since 2-good rings are closed under factor rings. By [12, Proposition 3], $RC_n$ is also 2-good. If $n = mp^k$ where $char\overline{R} = p \neq 2$, $k \geq 1$, and $(m,\ p) = 1$. Then $C_n \cong C_{p^k} \times C_m$, and so $RC_n \cong (RC_{p^k})C_m$. By [9, Theorem], $RC_{p^k}$ is also a local ring and $charRC_{p^k} = p$. The rest is proved similarly as above since $(p,\ m) = 1$. Thus we complete the proof. $\square$

By Theorem 10, we obtain the following corollary immediately

**Corollary 11.** *Let $R$ be a local ring with $\overline{R} = R/J(R)$ and let $C_n$ be a cyclic group of order $n$. If $char\overline{R} \neq 2$, then $RC_n$ is 2-clean.*

**Corollary 12.** ([16, Theorem 2.3]) *If $C_3$ is a cyclic group of order 3, then the group ring $\mathbb{Z}_{(p)}C_3$ is 2-clean for any prime number $p \neq 2$.*

**Remark 13.** The group ring $RC_n$ which satisfies the conditions of Theorem 10 need not be clean. In [5], Han and Nicholson showed that the group ring $\mathbb{Z}_{(7)}C_3$ is not clean where $\mathbb{Z}_{(7)} = \{m/n \in \mathbb{Q} : 7 \nmid n\}$.

Let $C_m = \{1, g, g^2, \cdots, g^{m-1}\}$ with $g^m = 1$ where $m$ is odd. Set $S = \{1, 2, \cdots, m-1\}$. Define $\sigma : S \longrightarrow S$ by $i \longmapsto 2i \pmod{m}$. It is checked easily that $\sigma$ is a permutation of $\{1, 2, \cdots, m-1\}$. Let $F$ be a field with $charF = 2$ and let $e = e_0 + e_1 g + \cdots + e_{m-1} g^{m-1} \in FC_m$ be an idempotent. Note that $2 = 0$ and $g^n = 1$, so $e^2 = e_0^2 + e_{\sigma(1)} g^{\sigma(1)} + \cdots + e_{\sigma(m-1)} g^{\sigma(m-1)}$. Suppose that $\sigma$ is a cyclic permutation. Then we have $e_0^2 = e_0$ and $e_1^2 = e_1 = e_2 = \cdots = e_{m-1}$, and so idempotents of $FC_m$ are 0, 1, $1 + g + \cdots + g^{m-1}$, $g + g^2 + \cdots + g^{m-1}$.

**Theorem 14.** *Let $R$ be a local ring with $char\overline{R} = 2$ and let $C_n$ be a cyclic group of order $n$. Write $n = m \cdot 2^k$ ($k \geq 0$) where $(m, 2) = 1$. If $\overline{R}$ ia a field and $\sigma$ is a cyclic permutation of $\{1, 2, \cdots, m - 1\}$, then the group ring $RC_n$ is semiperfect.*

**Proof.** Suppose $k \geq 1$. Then $C_n \cong C_{2^k} \times C_m$ from $(m, 2) = 1$, and so $RC_n \cong (RC_{2^k})C_m$. By [9, Theorem], $RC_{2^k}$ is local. Since $\overline{R}$ is a field and $\overline{R}C_{2^k} \twoheadrightarrow \overline{RC_{2^k}}$ is a ring epimorphism, $\overline{RC_{2^k}}$ is a field and $char\overline{RC_{2^k}} = char\overline{R} = 2$. Hence we may assume $n = m$. Note that $\overline{R}C_m$ is semisimple by $(m, 2) = 1$ and $J(R)C_m \subseteq J(RC_m)$, so $J(R)C_m = J(RC_m)$. This shows that $\overline{RC_m} \cong \overline{R}C_m$ with $char\overline{R} = 2$. Since $\overline{R}$ is a field and $\sigma$ is a cyclic permutation of $\{1, 2, \cdots, m - 1\}$, $\overline{R}C_m$ has only four idempotents, and so all idempotents in $\overline{RC_m}$ are $\overline{0}$, $\overline{1}$, $\overline{1} + \overline{g} + \cdots + \overline{g}^{m-1}$, $\overline{g} + \overline{g}^2 + \cdots + \overline{g}^{m-1}$. We find easily idempotents in $RC_m$, $f_1 = 0$, $f_2 = 1$, $f_3 = m^{-1}(1 + g + \cdots + g^{m-1})$, $f_4 = m^{-1}((m - 1) - g - g^2 - \cdots - g^{m-1})$ such that $\overline{f}_1 = \overline{0}$, $\overline{f}_2 = \overline{1}$, $\overline{f}_3 = \overline{1} + \overline{g} + \cdots + \overline{g}^{m-1}$, $\overline{f}_4 = \overline{g} + \overline{g}^2 + \cdots + \overline{g}^{m-1}$. This shows that $RC_m$ is semiperfect. □

The following result is immediate by Theorem 14 and by [1, Theorem 9].

**Corollary 15.** *Let $R$ be a local ring with $char\overline{R} = 2$ and let $C_n$ be a cyclic group of order $n$. Write $n = m \cdot 2^k$ ($k \geq 0$) where $(m, 2) = 1$. If $\overline{R}$ is a field and $\sigma$ is a cyclic permutation of $\{1, 2, \cdots, m - 1\}$, then the group ring $RC_n$ is clean.*

**Corollary 16**. ([16, Theorem 3.2]) *If $C_3$ is a cyclic group of order 3, then the group ring $\mathbb{Z}_{(2)}C_3$ is clean.*

**Remark 17.** The condition which $\sigma$ is cyclic in Theorem 14 can not be removed. In fact, it is determined only by $m$ whether the permutation $\sigma$ of $\{1, 2, \cdots, m - 1\}$ is cyclic. We calculate that $\sigma$ is cyclic in the case $m = 3, 5, 11, 13, \cdots$. However, set $m = 7$ or 9, $\sigma$ is not cyclic. Here, $\mathbb{Z}_{(2)}C_7$ is not semiperfect. In fact, in $\mathbb{Z}_2[X]$, $X^7 - \overline{1} = (X + \overline{1})(X^3 + X - \overline{1})(X^3 + X^2 + \overline{1})$. But in $\mathbb{Z}_{(2)}[X]$, $X^7 - 1 = (X - 1)(X^6 + X^5 + X^4 + X^3 + X^2 + X + 1)$ and $X^6 + X^5 + X^4 + X^3 + X^2 + X + 1$ is irreducible. So $\mathbb{Z}_{(2)}C_7$ is not semiperfect by [15, Theorem 5.8]. Note that $\overline{\mathbb{Z}_{(2)}C_7}$ is semisimple, then idempotents cannot be lifted modulo $J(\mathbb{Z}_{(2)}C_7)$, and so $\mathbb{Z}_{(2)}C_7$ is not clean.

*Department of Mathematics*
*Southeast University*
*Nanjing, 210096, China*
*e-mail: fylwangz@126.com*
*jichen@seu.edu.cn*